\documentclass[12pt,reqno]{amsart}
\usepackage{amsmath}
\usepackage{amssymb}
\usepackage[all]{xy}

\newif\ifprivate
\privatefalse

\def\???{\ifprivate {\bf {???}} \marginpar{{\Huge {\bf ?}}}
\else \fi}

\def\???{{\bf {???}} \marginpar{{\Huge {\bf ?}}} }
\newcount\minleft
\newcount\timehour
\def\thetime{\timehour=\time
\divide\timehour by60 \minleft=\timehour \multiply\minleft by -60
\advance\minleft by\time \ifnum\time>720\advance\timehour by-12\fi\relax
\number\timehour:\ifnum\minleft<10 %
    0\fi\relax\number\minleft
    \ifnum\time>720~pm \else~am\fi}

\newtheorem{defn0}{Definition}[section]
\newtheorem{prop0}[defn0]{Proposition}
\newtheorem{thm0}[defn0]{Theorem}
\newtheorem{lemma0}[defn0]{Lemma}
\newtheorem{corollary0}[defn0]{Corollary}
\newtheorem{example0}[defn0]{Example}
\newtheorem{remark0}[defn0]{Remark}
\newtheorem{conjecture0}[defn0]{Conjecture}

\newenvironment{definition}{ \begin{defn0}}{\end{defn0}}
\newenvironment{proposition}{\bigskip \begin{prop0}}{\end{prop0}}
\newenvironment{theorem}{\bigskip \begin{thm0}}{\end{thm0}}
\newenvironment{lemma}{\bigskip \begin{lemma0}}{\end{lemma0}}
\newenvironment{corollary}{\bigskip \begin{corollary0}}{\end{corollary0}}

\newenvironment{remark}{ \begin{remark0}\rm}{\end{remark0}}
\newenvironment{conjecture}{\begin{conjecture0}}{\end{conjecture0}}

\newcommand{\propref}[1]{Proposition~\ref{#1}}
\newcommand{\thmref}[1]{Theorem~\ref{#1}}
\newcommand{\lemref}[1]{Lemma~\ref{#1}}
\newcommand{\corref}[1]{Corollary~\ref{#1}}

\newcommand{\remref}[1]{Remark~\ref{#1}}

\def\max{{\bf m}}                   
\def\res{{\bf k}}                   
\newcommand{\MIJ}{\mathcal{M}^{I,J}}
\newcommand{\SIJ}{\Sigma^{I,J}}
\newcommand{\KIJ}{K^{I,J}}
\def\rees{{\mathcal R}}             

\begin{document}
\begin{flushleft}{\textsf{\tiny{To apperar in The Royal Society of Edinburgh Proceedings (Section A) Mathematics}}}\end{flushleft}

\title[Bigraded structures]{{\bf
Bigraded structures and the depth of blow-up algebras}}
\author[Gemma Colom\'e-Nin]{Gemma Colom\'e-Nin ${}^{*}$}
\thanks{${}^{*}$
FPI grant  of DGICYT BFM2001-3584}
\author[Juan Elias]{Juan Elias ${}^{**}$}
\thanks{${}^{*}$ ${}^{**}$Partially supported by MEC-FEDER MTM2004-01850.\\
\rm \indent 2000 MSC:  13A30, 13C14,  13D40}
\address{Departament d'\`{A}lgebra i Geometria
\newline \indent Facultat de Matem\`{a}tiques
\newline \indent Universitat de Barcelona
\newline \indent Gran Via 585, 08007
Barcelona, Spain} \email{{\tt gcolome@ub.edu}} \email{{\tt elias@ub.edu}}
\date{November 17, 2005}

\begin{abstract}
Let $R$ be a Cohen-Macaulay local ring, and let $I\subset R$ be an ideal with
minimal reduction $J$. In this paper we attach to the pair $I$, $J$ a
non-standard bigraded module $\Sigma^{I,J}$.
 The study of the bigraded
Hilbert function of $\SIJ$ allows us to prove  a improved  version of Wang's
conjecture  and a weak version of Sally's conjecture, both on the depth of the
associated graded ring $gr_I(R)$. The module $\SIJ$  can be considered as a
refinement of the Sally's module previously introduced by W. Vasconcelos.
\end{abstract}

\maketitle

\baselineskip 16pt

\section*{Introduction}

 Let  $(R,\max,\res) $ be a $d$-dimensional
Cohen-Macaulay local ring. Let $I$ be an $\max-$primary ideal of $R$ with
minimal reduction $J$. One of the major problems in commutative algebra is to
estimate the depth of the associated graded ring $gr_I(R)=\bigoplus_{n\ge 0}
I^{n}/I^{n+1}$ and the Rees algebra $\rees(I)=\bigoplus_{n\ge 0} I^n t^n$ for
ideals $I$ having good properties. Attached to the pair $I,J$ we can consider
the integers
$$
\Delta(I,J)=\sum_{p\ge 1}length_R\left(\frac{I^{p+1}\cap J}{I^{p}J}\right)
\text{\; \quad , \quad }
 \Lambda(I,J)=\sum_{p\ge 0}length_R\left(\frac{I^{p+1} }{J I^{p}}\right),
 $$
$\Delta_p(I,J)=length_R(I^{p+1}\cap J/J I^p)$, and
$\Lambda_p(I,J)=length_R(I^{p+1}/J I^p)$, for $p\ge 0$.

Related to these integers there are some results and conjectures on the depth
of $gr_I(R)$ that we next review. Valabrega and Valla proved that if
$\Delta(I,J)=0$ then $gr_I(R)$ is Cohen-Macaulay, \cite{VV78}. Based in this
result Guerrieri proposed the following conjecture:

 \begin{conjecture}[Guerrieri, \cite{Gue94}]
Let $I$ be an $\max-$primary ideal of $R$ with minimal reduction $J$. Then
$$
depth(gr_I(R)) \ge d-\Delta(I,J).
$$
\end{conjecture}
\noindent
 Guerrieri proved
the case $\Delta(I,J)=1$ and some partial cases for $\Delta(I,J)=2$,
\cite{Gue94}. Wang proved the case $\Delta(I,J)=2$ without any restriction,
\cite{Wan00}.

Guerrieri in her thesis asked if the conditions $\Delta_p(I,J)\le 1$, $p\ge 1$,
implies that $depth(gr_I(R)) \ge d-1$, \cite{Gue93}, Question 2.23. Wang in
\cite{Wan02b}, Example 3.13, gave a counterexample to Guerrieri's question and
asked if this question has an affirmative answer assuming that $R$ is a regular
local ring.

\medskip
Huckaba and Marley proved that $e_1(I)\le \Lambda(I,J)$ and if the equality
holds then $depth(gr_I(R))\ge d-1$, \cite{HM95}. Hence we can consider the
non-negative integer $\delta(I,J)=\Lambda(I,J)-e_1(I)\ge 0$. Wang  showed that
$\delta(I,J) \le \Delta(I,J)$ and that Guerrieri's conjecture is implied by the
following one, \cite{Wan00},
\begin{conjecture}[Wang]
Let $I$ be an $\max-$primary ideal of $R$ with minimal reduction $J$. Then
$$
depth(gr_I(R)) \ge d-1-\delta(I,J).
$$
\end{conjecture}
\noindent
 Huckaba proved
the conjecture in the case $\delta(I,J)=0$, \cite{Huc96}, \cite{HM95}.
  If
$\delta(I,J)=1$ Wang proved the conjecture and Polini gave a simpler proof,
\cite{Wan00}, \cite{Pol00}. For $\delta(I,J)=2$ Rossi and Guerrieri proved
Wang's conjecture assuming that $R/I$ is Gorenstein, \cite{GR99}. Wang gave a
counterexample to the conjecture for $d=6$,  \cite{Wan01}.

In the main result of this paper we prove a refined version of Wang's
conjecture, \thmref{generalitzacio}. We naturally decompose the integer
$\delta(I,J)=\sum_{p\ge0}\delta_p(I,J)$ as a finite sum of non-negative
integers $\delta_p(I,J)$, with $\Delta_p(I,J)\ge \delta_p(I,J)\ge 0$, see
section two. Let us consider
 the maximum, say $\bar{\delta}(I,J)$, of the integers $\delta_p(I,J)$ for
$p\ge 0$.

\bigskip
\noindent {\bf \thmref{generalitzacio}.}{ \em Assume that $\bar{\delta}(I,J)\le
1$. Then
$$
depth(\rees(I))\ge d- \bar{\delta}(I,J)
$$
and $depth(gr_I(R))\ge d-1-\bar{\delta}(I,J).$ }

\bigskip
\noindent There is another conjecture that considers some conditions on the
modules $I^{p+1}/JI^p$ and $I^p\cap J/I^{p-1}J$, it is Sally's conjecture:

\begin{conjecture}[Sally]
Let $I$ be an $\max-$primary ideal of $R$ with minimal reduction $J$. If $I^n
\cap J= I^{n-1}J$ for $n=2,\ldots,t$ and $length(\frac{I^{t+1} }{J
I^{t}})=\epsilon \le \text{ Min} \{1, d-1\}$ then it holds
$$
d-\epsilon \le depth(gr_I(R)) \le d.
$$
\end{conjecture}
\noindent This conjecture was proved by Corso-Polini-Vaz-Pinto,
 Elias,  and Rossi, \cite{CPP97}, \cite{Eli99}, \cite{Ros00}.

\bigskip
The aim of this work is to introduce  a non-standard bigraded module $\SIJ$ in
order to study the depth of the associated graded ring $gr_I(R)$ and the Rees
algebra $\rees(I)$ of $I$. A secondary purpose is to present a unified
framework where several results and objects appearing in the papers on the
above conjectures  can be studied,  \remref{unificacio}. The key tool of this
paper is the Hilbert function of non-standard bigraded modules.

The first section is mainly devoted to recall some preliminary results on Sally
module and the Hilbert function of non-standard bigraded modules.

In Section two we introduce a non-standard bigraded module $\Sigma^{I,J}$
naturally attached  to $I$ and a minimal reduction $J$ of $I$, this module can
be considered as a refinement of the Sally module previously introduced by W.
Vasconcelos. From a natural presentation of $\SIJ$ we define two  bigraded
modules $\KIJ$ and  $\MIJ$,  and we consider some diagonal submodules of
 them: $\SIJ_{[p]}$ and $\KIJ_{[p]}$.

For all $p\ge 0$ we consider the integer $\delta_p(I,J)=e_0(\KIJ_{[p]})$. In
\propref{localversion} we prove following inequalities
$$
\Delta_p(I,J) \ge \delta_p(I,J)=\Lambda_p(I,J)-e_0(\SIJ_{[p]})\ge 0.
$$
\noindent Summing up these inequalities with respect $p$ we recover the
inequalities, \cite{Wan00}, \cite{HM95},
$$
\Delta(I,J) \ge \delta(I,J)=\Lambda(I,J)-e_1(I)\ge 0.
$$
In particular we decompose the first Hilbert coefficient $e_1(I)$ of $I$ as a
sum of the multiplicities $e_0(\SIJ_{[p]})$ when $p$ ranges the set of
non-negative integers:
$$
e_1(I)=\sum_{p\ge 0} e_0(\SIJ_{[p]}).
$$

Section three is devoted to prove a refined version of  Wang's conjecture by
considering some special configurations of the set $\{\delta_p(I,J)\}_{p\ge 0}$
instead of $\delta(I,J)=\sum_{p\ge 0} \delta_p(I,J)$,  \thmref{generalitzacio}.
This version allows us to improve the bound for the $depth$  in Guerrieri's
question, \propref{GUE-question} (see \cite{Gue93}, Question 2.23). As a
by-product we recover the known cases of Wang's conjecture, \corref{wang}, and
we prove a weak version of Sally's conjecture, \corref{sallyconj}. An essential
point of this section is based on generalize a part of the work of Polini in
\cite{Pol00} following her ideas.

\medskip
\noindent {\sc acknowledgements.} The authors thank M.E. Rossi for the
discussions regarding the main ideas of this paper.


\newpage
\bigskip
\section{Preliminaries}

Let $(R,\max)$ be a local ring of dimension $d>0$, and let $I$ be an
$\max$-primary ideal of $R$. Without loss of generality we may assume that the
residue field $\res = R/ \max $ is infinite.
 We set $I^i=0$ for $i<0$, and  $I^i=R$ for $i=0$.

The  Rees algebra of $R$ associated to $I$ is the $R-$algebra
$$\rees(I)=\bigoplus_{n\ge 0} I^n t^n$$
and the associated graded ring of $R$ with respect to $I$ is
$$gr_I(R)=\bigoplus_{n\ge 0}\frac{I^n}{I^{n+1}} t^n.$$

An ideal $J\subseteq I$ is said to be a reduction of $I$ if there exists an
integer $r\ge 0$ such that $I^{r+1}=JI^r$. $J$ is a minimal reduction of $I$ if
$J$ is a reduction of $I$ an $J$ itself does not contain any proper reduction.
If $J$ is a minimal reduction of $I$, the  reduction number of $I$ with respect
to $J$ is the least integer $r_J(I)$ such that $ I^{r+1}=JI^r$ for all $r \ge
r_J(I)$.
 The reduction number  $r(I)$  of $I$ is defined as
lowest integer $r_J(I)$ where $J$ is a minimal reduction of $I$.

We denote by $h^0_{I}(n)=length_R(I^n/I^{n+1})$ the $0-th$ Hilbert function of
$I$. The higher Hilbert functions of $I$ are defined by, $i\ge 0$,
$$
h^{i+1}_{I}(n)=\sum_{j=0}^{n}h^i_{I}(j).
$$

\noindent It is well known that there exist integers $e_j(I) \in \mathbb Z$
such that if we write
$$
p^i_{I}(X)=\sum_{j=0}^{d+i-1} (-1)^j e_j(I) \binom{X+d+i-j-1}{d+i-j-1}
$$

\noindent then $p^i_{I}$ is the $i-$th Hilbert polynomial of $I$, i.e.
$h^i_{I}(n)=p^i_{I}(n)$ for $n \gg 0$. The integer $e_i(I)$ is the $i-$th
Hilbert coefficient of $I$. In this paper we set $h_I=h_I^0$ and $p_I=p^0_I$.
We denote by $pn(I)$ the postulation number of $h_{I}$, i.e. the least integer
$t$ such that $h_{I}(t+n)=p_{I}(t+n)$ for all $n \ge 0$.

 Let $I$ an
$\max$-primary ideal of $R$ and $J$ a minimal reduction of $I$. The  Sally
module of $I$ with respect to $J$ is the $\rees(J)$-module
$$
S_J(I)=\frac{I\rees(I)}{I\rees(J)}=\bigoplus_{n\ge 1}\frac{I^{n+1}}{J^nI} t^n
$$

\noindent
 We define the Hilbert
function of the Sally module $S_J(I)$ as, \cite{Vas94},
$$h_{S_J(I)}(n)=length_R({I^{n+1}}/{J^nI}).$$
If $S_J(I)\ne 0$ then  $dim(S_J(I))=d$, and we can consider the Hilbert
polynomial of $S_J(I)$
$$p_{S_J(I)}(n)=\sum_{i=0}^{d-1}(-1)^is_i\binom{n+d-i-1}{d-i-1}.$$

In the next proposition we collect some known results on the Sally module.

\begin{proposition}
Let $(R,\max)$ a Cohen-Macaulay local ring of dimension $d>0$. Let $I$ be an
$\max$-primary ideal of $R$ and $J$ a minimal reduction of $I$.
\begin{enumerate}
\item[(i)] If $S_J(I)=0$ then $gr_I(R)$ is Cohen-Macaulay. If $S_J(I)=0$ and
$d\ge 2$, then $\rees(I)$ is Cohen-Macaulay. \item[(ii)] $depth(gr_I(R))\ge
depth(S_J(I))-1$. \item[(iii)] If $depth(gr_I(R))<d$, then
$depth(S_J(I))=depth(gr_I(R))+1$. \item[(iv)] The Hilbert coefficients of
$S_J(I)$ and $I$ are related by the following equalities

$$\left\{\begin{array}{ll} e_0(I)=length_R(R/J) & \\ \\
e_1(I)=e_0(I)-length_R(R/I)+s_0=s_0+length_R(I/J) & \\ \\
e_i(I)=s_{i-1} & i=2,\ldots,d \end{array}\right.$$
\end{enumerate}
\end{proposition}

\bigskip
 Roberts proved the existence
of Hilbert polynomials of bigraded modules over bigraded rings
$\res[X_1,\dots,X_s]$ with  variables $X_1,\dots,X_s$ of bidegrees $(1,0)$,
$(0,1)$, and  $(1,1)$, here $\res$ is a field, \cite{Rob00}, see also
\cite{HoTr02}. The results of Section 3 of \cite{Rob00} can be easily
generalized to polynomial rings with coefficients in an Artin ring.

Let $A=C[X_0,\dots,X_s,T_0,\dots,T_u,S_1,\dots,S_v]$ be a bigraded polynomial
ring  over an Artin ring $C$ in the variables $X_0,\dots,X_s,T_0,\dots,T_u$ and
$S_1,\dots,S_v$. We assume that the variables $X_i$ have bidegree $(1,0)$, the
variables $T_i$ have bidegree $(1,1)$, and the variables  $S_i$ have  bidegree
$(0,1)$.

For a  bigraded $A$-module $M$ and for any $m,n\in\mathbb{Z}$, let $M_{(m,n)}$
be the piece of $M$ of bidegree $(m,n)$.
 Let $h_M(m,n)$ be the Hilbert function of $M$
defined by the equality
$$h_M(m,n)=\sum_{i\le n}length_A(M_{(m,i)}).$$

\begin{theorem}\label{polin bigrad}
Let $A=C[X_0,\dots,X_s,T_0,\dots,T_u,S_1,\dots,S_v]$ be a bigraded polynomial
ring  over an Artin ring $C$ in variables $X_0,\dots,X_s,T_0,\dots,T_u$ and
$S_1,\dots,S_v$, where each $X_i$ has bidegree $(1,0)$, each $T_i$ has bidegree
$(1,1)$, and each $S_i$ has bidegree $(0,1)$. For all   finitely generated
bigraded $A$-module $M$,
 there exist a polynomial in two
variables $p_M(m,n)$, and integers $m_0$ and $n_0$  such that
$$p_M(m,n)=h_M(m,n)$$ for all $(m,n)$ with $m\ge m_0$ and $n\ge m+n_0$.
\end{theorem}

\bigskip
\section{Bigraded Sally module}

Let $R$ be a Cohen-Macaulay local ring. Let $I=(b_1,\ldots,b_{\mu})$ be an
$\max$-primary ideal of $R$ and let $J=(a_1,\dots,a_d)$ be a minimal reduction
of $I$. Since $Jt \rees(I)$ is an homogeneous ideal of the  graded ring
$\rees(I)$ we can consider the associated graded ring of $\rees(I)$ with
respect to the homogeneous ideal $Jt \rees(I)=\bigoplus_{n
  \ge 0} J I^{n-1} t^n$
$$
gr_{Jt}(\rees(I))= \bigoplus_{j\ge 0} \frac{(Jt \rees(I))^{j}}{(Jt
\rees(I))^{j+1}} U^j.
$$
This ring has a natural bigraded structure that we describe briefly. Notice
that $\frac{\rees(I)}{Jt\rees(I)}=\bigoplus_{i\ge 0}\frac{I^i}{I^{i-1} J} t^i $
is a homomorphic image of the graded ring $R[V_1,\dots, V_{\mu}]$
 by the degree one  $R$-algebra homogeneous morphism
$$
\sigma: R[V_1,\dots, V_{\mu}] \longrightarrow
\frac{\rees(I)}{Jt\rees(I)}=\bigoplus_{i\ge 0}\frac{I^i}{I^{i-1} J} t^i
$$
\noindent defined by $\sigma(V_i)=b_i t \in \frac{I}{J} t$; $R[V_1,\dots,
V_{\mu}]$ is endowed with the standard graduation.

Let us consider the bigraded ring $B:=R[V_1,\dots, V_{\mu};T_1,\dots,T_d]$ with
$deg(V_i)=(1,0)$ and $deg(T_i)=(1,1)$, then there exists an exact sequence of
bigraded $B$-rings
\begin{equation}
\label{fond} 0 \longrightarrow K^{I,J} \longrightarrow
C^{I,J}:=\frac{\rees(I)}{Jt\rees(I)}[T_1,\dots,T_d]
\stackrel{\pi}{\longrightarrow} gr_{Jt}(\rees(I))\longrightarrow 0
\end{equation}

\noindent with $\pi(T_i)=a_i t U$, $i=1,\dots,d$; $K^{I,J}$ is the ideal of
initial forms of $Jt \rees(I)$. The $(i+j, j)$-graded piece of
$gr_{Jt}(\rees(I))$ and  $ C^{I,J}$ are
$$
gr_{Jt}(\rees(I))_{(i+j,j)}= \frac{I^i J^j }{ I^{i-1} J^{j+1}} t^{i+j} U^j
\quad \text{and} \qquad
C^{I,J}_{(i+j,j)}=\frac{I^i}{JI^{i-1}}t^{i}[T_1,\dots,T_d]_j,
$$
respectively. Notice that we have an $R$-algebra isomorphism
$\phi:R[T_1,\dots,T_d]/(\{a_iT_j-a_jT_i\}_{i,j}) \cong R[JtU]=\rees(J)$ defined
by $\phi(\overline{T_i})=a_i t U$, $i=1,\dots,d$. Observe that we write $tU$
instead of only $t$ to bear in mind the the bigraduation.

\bigskip
Given a $B-$bigraded module $M$ and an integer $p\in \mathbb Z$, we denote by
$M_{[p]}$
 the additive sub-group of $M$ defined
by the direct sum of the pieces $M_{(m,n)}$ such that $m-n=p+1$. Notice that
the product by the variable  $T_i$ induces an endomorphism of
$R[T_1,\dots,T_d]-$modules $M_{[p]} \stackrel{T_i}{\longrightarrow} M_{[p]}$,
and the product by $V_j$ a morphism of $R[T_1,\dots,T_d]-$modules $M_{[p]}
\stackrel{V_j}{\longrightarrow} M_{[p+1]}$. Hence $M_{\ge p}=\bigoplus_{n\ge
p}M_{[n]}$ is a sub-B-module of $M$, and we can consider  the exact sequence of
$R[T_1,\dots,T_d]$-modules
$$
0\longrightarrow M_{[p]} \longrightarrow M_{\ge p} \longrightarrow M_{\ge p+1}
\longrightarrow 0.$$

Moreover, in our case, the modules $\KIJ_{[p]}$, $C^{I,J}_{[p]}$ and
$gr_{Jt}(\rees(I))_{[p]}$ are $\rees(J)-$modules.

\bigskip
Next lemma shows that $\KIJ_{[p]}$, $C^{I,J}_{[p]}$ and
$gr_{Jt}(\rees(I))_{[p]}$ do not eventually vanish for a finite set of indexes
$p\in \mathbb Z$.

\begin{lemma}
\label{vanishing} $(i)$ For all $p\le -2$ or $p\ge r_J(I)$,  $C^{I,J}_{[p]}=0$,
$gr_{Jt}(\rees(I))_{[p]}=0$ and $K^{I,J}_{[p]}=0$.

\noindent $(ii)$ $\pi$ induces the following isomorphisms of
$\rees(J)-$modules:
$$gr_{Jt}(\rees(I))_{[0]}\cong
C^{I,J}_{[0]}\cong\frac{I}{J} t [T_1,\ldots,T_d],$$
$$gr_{Jt}(\rees(I))_{[-1]}\cong \rees(J),\quad C^{I,J}_{[-1]}\cong R[T_1,\dots,T_d].$$
Moreover, $K^{I,J}_{[0]}=0$.

\end{lemma}
\begin{proof}
In order to prove $(i)$, first we observe how  $C^{I,J}_{[p]}$ and
$gr_{Jt}(\rees(I))_{[p]}$ are. Since $C^{I,J}\cong \bigoplus_{i\geq
0}\frac{I^i}{JI^{i-1}}t^i[T_1,\dots,T_d]$, we have that
$$
C^{I,J}_{[p]}=\bigoplus_{m-n=p+1}C^{I,J}_{(m,n)}=\frac{I^{p+1}}{JI^p}t^{p+1}[T_1,\dots,T_d],
$$
and
$$
gr_{Jt}(\rees(I))_{[p]}=\bigoplus_{m-n=p+1}gr_{Jt}(\rees(I))_{(m,n)}=\bigoplus_{i\geq
0}\frac{J^iI^{p+1}}{J^{i+1}I^p}t^{p+1+i}U^i.
$$

Since $I^{i}=0$ for all $i<0$, we have that  $C^{I,J}_{[p]}=0$ and
$gr_{Jt}(\rees(I))_{[p]}=0$ for all $p\le -2$. By the definition of $r_J(I)$ we
have $I^{p+1}=JI^p$ for all $p\geq r_J(I)$, so
 $C^{I,J}_{[p]}=0$ and $gr_{Jt}(\rees(I))_{[p]}=0$ for
all $p\geq r_J(I)$. Notice that we have that $K^{I,J}_{[p]}\subseteq
C^{I,J}_{[p]}$ for each $p\in\mathbb{Z}$. Therefore, we have that
$K^{I,J}_{[p]}=0$ for all $p\leq -2$ and $p\geq r_J(I)$.

\medskip
\noindent $(ii)$  When $p=0$ we have
$$
C^{I,J}_{[0]}=\frac{I}{J}t[T_1,\dots,T_d]=It(R/J)[T_1,\dots,T_d]
$$
\noindent and
$$
gr_{Jt}(\rees(I))_{[0]}=\bigoplus_{i\geq 0}\frac{J^iI}{J^{i+1}}t^{i+1}U^i \cong
It gr_{J}(R).
$$
Since  $gr_J(R)\cong (R/J)[T_1,\dots,T_d]$, clearly we have that
$C^{I,J}_{[0]}\cong gr_{Jt}(\rees(I))_{[0]}$. By the exact sequence
\eqref{fond} we deduce  $K^{I,J}_{[0]}=0$. If $p=-1$ then we have
$C^{I,J}_{[-1]}=R[T_1,\dots,T_d]$ and
$gr_{Jt}(\rees(I))_{[-1]}=\bigoplus_{i\geq 0}J^it^iU^i \cong \rees{(J)}$.
\end{proof}

\bigskip
Let us consider the following bigraded finitely generated $B$-modules:

\vskip 2mm \noindent $\Sigma^{I,J}:=\bigoplus_{p\ge 0}
gr_{Jt}(\rees(I))_{[p]}$,

\vskip 2mm \noindent $\mathcal M^{I,J}:=\bigoplus_{p\ge 0} C^{I,J}_{[p]}\cong
\bigoplus_{p\ge 0} I^{p+1}/I^p J \ t^{p+1} [T_1,\ldots,T_d]$, \vskip 2mm
\noindent and from now on we consider the new
 \vskip 2mm
\noindent $\KIJ:=\bigoplus_{p\ge 0} K^{I,J}_{[p]}$.

\vskip 2mm \noindent Notice that  by \lemref{vanishing} there exists a natural
isomorphism of $\rees(J)-$modules
$$gr_{Jt}(\rees(I)) \cong \rees(J) \oplus \Sigma^{I,J}.$$

Since the modules $\SIJ$ and $\MIJ$ are annihilated by $J$, from the exact
sequence \eqref{fond} we deduce
 the  following exact sequence of
 $A=B\otimes_R R/J\cong R/J[V_1,\ldots,V_{\mu};T_1,\ldots,T_d]-$bigraded
modules
\begin{equation*}
\tag{S} 0 \longrightarrow K^{I,J} \longrightarrow \mathcal M^{I,J}
\longrightarrow \Sigma^{I,J}\longrightarrow 0.
\end{equation*}
\noindent
 Notice that from \lemref{vanishing} all relevant information
 of \eqref{fond}  is encoded in the exact sequence  (S).

For all $p\ge 0$ we have an exact sequence of $R/J[T_1,\ldots,T_d]-$modules
\begin{equation*}
\tag{$\textrm{S}_{[p]}$} 0 \longrightarrow \KIJ_{[p]} \longrightarrow
\MIJ_{[p]}=\frac{I^{p +1}}{J I^{p}}[T_1,\ldots,T_d] \longrightarrow \SIJ_{[p]}
\longrightarrow 0,
\end{equation*}
so we can consider the (classic) Hilbert function of $\SIJ_{[p]}$, $\MIJ_{[p]}$
and $\KIJ_{[p]}$ with respect the variables $T_1,\ldots,T_d$.

\begin{definition}
Let $I$ be an $\max$-primary ideal and let $J$ be a minimal reduction of $I$,
$\Sigma^{I,J}$ is the bigraded Sally module of I with respect $J$.
\end{definition}

\bigskip
In the next result we show some properties of the $A$-module $\MIJ$ that will
be used in this paper.

\begin{proposition}
\label{assoc} Given an $\max$-primary ideal $I$ of $R$  with minimal reduction
$J$. Then $Ass_A(\mathcal M^{I,J})= \{ \max A + (V_1,\ldots,V_{\mu})\}$ and
$\mathcal M^{I,J}$ is a Cohen-Macaulay $A-$module of dimension $d$.
\end{proposition}

\begin{proof}
By the definitions of $Ass_A(\MIJ)$ and $\MIJ$, with a simply computation we
can show that $Ass_A(\MIJ)=\{Q\}$, where $Q=\max A +(V_1,\ldots,V_{\mu})$.
Similarly, it is easy to see that $Ann_A(\MIJ)=Q$. Notice that $A/Q\cong
(R/\max)[T_1,\ldots,T_d]$. Hence,
$$dim
(\MIJ)=dim(A/Ann(\MIJ))=dim(A/Q)=dim((R/\max)[T_1,\dots,T_d])=d.$$ Since
$T_1,\dots,T_d$ is a regular sequence in $A$ we have that $\MIJ$ is a
Cohen-Macaulay module of dimension $d$.
\end{proof}

\bigskip
\begin{remark}
\label{link} Notice that the length as $R-$module of
$$
\Sigma^{I,J}_{(m+1,*)}\cong \bigoplus_{j=0}^{j=m}
\frac{I^{m+1-j}J^j}{I^{m-j}J^{j+1}}t^{m+1}U^{j}
$$
is equal to the length of
$$
S_J(I)_m \oplus \frac{I J^m}{J^{m+1}},
$$
where $S_J(I)_m$ is the degree $m$ piece of the Sally module $S_J(I)$.

 On the other hand
$$
\Sigma^{I,J}_{[p]}=gr_{Jt}(\rees(I))_{[p]}= \bigoplus_{i\ge
0}\frac{J^{i}I^{p+1}}{J^{i+1}I^{p}}t^{p+1+i}U^i
$$
 is isomorphic as
$\rees(J)-$module to the module $L_p$ defined by Vaz-Pinto in \cite{Vaz95} for
all $p\ge 1$.

In several papers of Wang appear the modules $T_{k,n}, S_{k,n}$ defined as the
kernel and co-kernel of some exact sequence, \cite{Wan97d}, \cite{Wan97a},
\cite{Wan00}, \cite{Wan01} and \cite{Wan02b}. One of the key results of Wang's
papers is to prove that there exists the Hilbert function of $T_{k,n}$ and
$S_{k,n}$. In our framework these results follows from the following
commutative diagram:
$$
\xymatrix{
 & & \bigoplus_{j=1}^{\binom{n+d-1}{d-1}}  \frac{I^j}{J I^{j-1}} \ar[dd]_{\cong} \ar[dr] & &\\
  0 \ar[r] & T_{k,n}= K^{I,J}_{(k+n,n)} \ar[dr] \ar[ur] & &
   S_{k,n}=\SIJ_{(k+n,n)}   \ar[r] & 0 \\
   & & \MIJ_{(k+n,n)} \ar[ur] &  & }
$$

Finally the $\rees(J)-$module $C_p=\bigoplus_{i\ge 1}\frac{I^{i+p}}{J^{i}
I^p}$,
 defined by Vaz-Pinto, can be linked to $\SIJ_{\ge p}$ by means of the following
 sets of exact sequences of $\rees(J)-$modules, \cite{Vaz95}.
 Here we set $r=r_J(I)$.
There exist two sets of exact sequences of $\rees(J)-$modules:

\begin{equation*}
\tag{SQ1}
\begin{array}{ccccccccc}
0 & \longrightarrow & \SIJ_{[1]} &\longrightarrow & S_J(I)=C_1 &
\longrightarrow &
C_2 & \longrightarrow & 0  \\ \\
0 & \longrightarrow & \SIJ_{[2]} &\longrightarrow & C_2 & \longrightarrow &
C_3 & \longrightarrow & 0  \\
    &               & \vdots &            & \vdots & &
\vdots & &   \\
0 & \longrightarrow & \SIJ_{[r-2]}& \longrightarrow & C_{r-2} & \longrightarrow
&
C_{r-1}=\SIJ_{[r-1]} & \longrightarrow & 0  \\
\end{array}
\end{equation*}

\bigskip
\noindent and

\begin{equation*}
\tag{SQ2}
\begin{array}{ccccccccc}
0 & \longrightarrow &\SIJ_{[0]}=0 &\longrightarrow & \SIJ= \SIJ_{\ge 0} &
\longrightarrow &
\SIJ_{\ge 1} & \longrightarrow & 0  \\ \\
0 & \longrightarrow & \SIJ_{[1]} &\longrightarrow & \SIJ_{\ge 1} &
\longrightarrow &
\SIJ_{\ge 2} & \longrightarrow & 0  \\
    &               & \vdots&             & \vdots & &
\vdots & &   \\
0 & \longrightarrow & \SIJ_{[r-2]}& \longrightarrow & \SIJ_{\ge r-2} &
\longrightarrow &
\SIJ_{\ge r-1}=\SIJ_{[r-1]} & \longrightarrow & 0  \\
\end{array}
\end{equation*}
\end{remark}


\bigskip

Now, given a bigraded $A-$module $M$,
$A=R/J[V_1,\ldots,V_{\mu};T_1,\ldots,T_d]$,
 we can consider the
Hilbert function of $M$ defined by
$$
h_M(m,n)=\sum_{0\le j \le n} length_A(M_{(m,j)})
$$
\noindent and by \thmref{polin bigrad}, there exist integers $f_{i,j}(M) \in
\mathbb Z$, $i\ge 0$, $j\ge 0$, and $i+j \le c-1$, for some integer $c\ge 0$,
such that the polynomial
$$
p_M(m,n)=\sum_{i+j\le c-1} f_{i,j}(M) \binom{m}{i} \binom{n}{j}
$$
\noindent verifies $p_M(m,n)=h_M(m,n)$ for all $m \ge m_0$ and $n \ge n_0+m$
for some integers $m_0, n_0 \ge 0$.

\bigskip
\begin{lemma}
\label{finitediag} Let $M$ be a bigraded $A-$module for which there exist
integers $a \le b$ such that
 $M_{[p]}=0$ for all $p \not\in [a,b]$,  then
 $f_{i,j}(M)=0$ for all $j\ge 1$.
\end{lemma}
\begin{proof}
Since  $M_{[p]}=\bigoplus_{m-n=p+1} M_{(m,n)}$. If $M_{[p]}=0$ for all $p\notin
[a,b]$, we have that $M_{(m,n)}=0$ for all pair $(m,n)$ such that $m-n>b+1$ or
$m-n<a+1$. If we take $m\ge m_0$, $n\ge n_0+m$ we have
$$h_M(m,n)=p_M(m,n)=\sum_{i+j\leq
c-1}f_{i,j}(M)\binom{m}{i}\binom{n}{j}$$

We can suppose that $m \ge b+1$ and $n\geq m-b-1$, so we have
\begin{eqnarray*}
h_M(m,n)&=&\sum_{0\leq j\leq n}length_A(M_{(m,j)})\\&=&\sum_{m-b-1 \leq j\leq
m-a-1}length_A(M_{(m,j)})=\sum_{j\in\mathbb{Z}}length_A(M_{(m,j)})
\end{eqnarray*}
because $M_{(m,j)}=0$ when $j<m-b-1$ and $j>m-a-1$. Hence, for $m,n$ large
enough we have that
$$h_M(m,n)=p_M(m,n)=h_M(m)=p_M(m).$$
Therefore
$$\sum_{i+j\leq
c-1}f_{i,j}(M)\binom{m}{i}\binom{n}{j}=\sum_i{a_i\binom{m}{i}}.
$$
Since $\{\binom{m}{i}\binom{n}{j}\}_{i,j}$ is a basis of the polynomial ring in
$m$ and $n$, we have that $f_{i,j}(M)=0$ for all $j\geq 1$.
\end{proof}

\bigskip

Given an $A-$module under the hypothesis of the above Lemma we can write
$$
p_M(m,n)=p_M(m)=\sum_{i=0}^{c-1} f_{i,0}(M) \binom{m+c-i}{c-i},
$$
 and
 $h_M(m,n)=p_M(m)$ for $m \ge m_0$ and $n \ge n_0+m$.
 From \lemref{vanishing} we can apply the last result to the
$A-$modules $\SIJ$, $\MIJ$, and $\KIJ$.

\begin{proposition}
\label{sigmahf} Let $I$ be an $\max-$primary ideal of $R$ with minimal
reduction $J$. Then
$$
p_{\SIJ}(m) = \sum_{i=0}^{d-1}(-1)^i e_{i+1}(I)\binom{m-1+d-i-1}{d-i-1}
$$
\end{proposition}
\begin{proof}
Since the length of $\Sigma^{I,J}_{(m+1,*)}$ is equal to the length of
$S_J(I)_m\oplus \frac{IJ^m}{J^{m+1}}$ as $R-$modules, we
 have that
$$
h_{\Sigma^{I,J}}(m,n)=length_R(S_{m-1})+ length_R(\frac{I
  J^{m-1}}{J^{m}})$$
for all $m\ge m_0$, $n\ge n_0+m$. Since $gr_J(R)\cong (R/J)[T_1,\dots,T_d]$,
clearly we get that $I gr_J(R)\cong (I/J)[T_1,\dots,T_d]$. Thus,
$length_R(\frac{IJ^{m-1}}{J^m})$ coincides with the length of the piece of
degree $m-1$ of  $(I/J)[T_1,\dots,T_d]$. So we have
$$
length_R(\frac{IJ^{m-1}}{J^m})=length_R(I/J)\binom{m-1+d-1}{d-1}.
$$
Hence we deduce that, \lemref{finitediag},
$$p_{\Sigma^{I,J}}(m)=p_{S_J(I)}(m-1)+
length_R(I/J)\binom{m-1+d-1}{d-1}.
$$
Let us recall that there exist integers $s_0, \ldots, s_{d-1}$ such that
$$
p_{S_J(I)}(n)= \sum_{i=0}^{d-1} (-1)^i s_i \binom{n+d-i-1}{d-i-1}
$$
and that $s_0=e_1(I)-length_R(I/J)$, $s_i=e_{i+1}(I)$ for $i=1,\ldots, d-1$,
\cite{Vas94}.
 Hence, we have
$$p_{S_J(I)}(m-1)=\sum_{i=0}^{d-1}(-1)^i
e_{i+1}(I)\binom{m-1+d-i-1}{d-i-1} - length_R(I/J)\binom{m-1+d-1}{d-1}$$ and
then
\begin{eqnarray} p_{\SIJ}(m) & = &
p_{S_J(I)}(m-1)+length_R(I/J)\binom{m-1+d-1}{d-1} \nonumber \\ & = &
\sum_{i=0}^{d-1}(-1)^i e_{i+1}(I)\binom{m-1+d-i-1}{d-i-1} \nonumber
\end{eqnarray}
\end{proof}

In the next proposition we compute the multiplicities of the modules $\MIJ$,
$\SIJ$ and $\KIJ$. Notice that they are related with the integer $\Lambda(I,J)$
defined in the introduction. From now on we consider the integer
$\delta(I,J)=\Lambda(I,J)-e_1(I)$.

\begin{proposition}
\label{multiplicities} The following conditions hold:

\noindent $(i)$ $deg(p_{\MIJ})=d-1$ and $e_0(\MIJ)= \Lambda(I,J)$.

\noindent $(ii)$ If $\SIJ =0$ then $gr_I(R)$ is a Cohen-Macaulay ring . If
$\SIJ \neq 0$ then $deg(p_{\SIJ})=d-1$ and $e_0(\SIJ)=e_1(I)$.

\noindent $(iii)$ $e_0(\KIJ)=\delta(I,J)$; if $\KIJ \neq 0$ then
$deg(p_{\KIJ})=d-1$. In particular,
$$\Lambda(I,J)\ge
e_1(I).$$
\end{proposition}
\begin{proof}
$(i)$ We know that  $\MIJ$ is Cohen-Macaulay of dimension $d$, \propref{assoc}.
By \lemref{vanishing} and \lemref{finitediag} we have that
$p_{\MIJ}(m,n)=p_{\MIJ}(m)$, so $deg(p_{\MIJ})=d-1$. Since  $\MIJ\cong
\bigoplus_{p\ge 0}I^{p+1}/I^pJ\ t^{p+1}[T_1,\dots,T_d]$ we have that
\begin{eqnarray*}
h_{\MIJ}(m,n) & = &
\sum_{i=0}^{n}length_R\left(\MIJ_{(m,i)}\right)=\sum_{i=0}^n
length_R\left(\frac{I^{m-i}}{I^{m-i-1}J}[T_1,\dots,T_d]_i\right) \\
& = & length_R\left(\frac{I^m}{I^{m-1}J}[T_1,\dots,T_d]_0\right)+
length_R\left(\frac{I^{m-1}}{I^{m-2}J}[T_1,\dots,T_d]_1\right)
\\ & & +\dots +
length_R\left(\frac{I^{m-n}}{I^{m-n-1}J}[T_1,\dots,T_d]_n\right)\end{eqnarray*}

\noindent Notice that for  $n\ge m \gg 0$ we have
\begin{eqnarray*}
p_{\MIJ}(m,n)=h_{\MIJ}(m,n) & = &
length_R\left(\frac{I^m}{I^{m-1}J}[T_1,\dots,T_d]_0\right)+ \\ & & +
length_R\left(\frac{I^{m-1}}{I^{m-2}J}[T_1,\dots,T_d]_1\right)+
\\ & & +\dots +
length_R\left(\frac{I}{J}[T_1,\dots,T_d]_{m-1}\right)\end{eqnarray*}

\noindent Moreover, for  $m\ge r_J(I)$ we have
\begin{eqnarray*}
p_{\MIJ}(m,n) & = & \sum_{i\ge
1}length_R\left(\frac{I^i}{I^{i-1}J}[T_1,\dots,T_d]_{m-i}\right)
\\ & = & \sum_{i\ge
1}length_R\left(\frac{I^i}{I^{i-1}J}\right)\binom{m-i+d-1}{d-1}\end{eqnarray*}
Each binomial number is a polynomial in $m$ of degree $d-1$, and the leading
coefficient of this polynomial $p_{\MIJ}$ gives us the multiplicity
$$e_0(\MIJ)=\sum_{i\ge 1}length_R\left(\frac{I^i}{I^{i-1}J}\right)=\Lambda(I,J)$$

\noindent $(ii)$ If $\SIJ =0$ then $S_J(I)=0$ and $gr_I(R)$ is Cohen-Macaulay,
\cite{Vaz95} \remref{link}. Let us assume that  $\SIJ \neq 0$. Since,
\propref{sigmahf},
$$
p_{\Sigma^{I,J}}(m)=\sum_{i=0}^{d-1}(-1)^i
e_{i+1}(I)\binom{m-1+d-i-1}{d-i-1},$$ then $deg(p_{\SIJ})=d-1$ and
$e_0(\SIJ)=e_1(I)$.

\noindent $(iii)$ From the additivity of the multiplicity in (S) and by $(i)$
and $(ii)$ we get
$$
e_0(\KIJ)=e_0(\MIJ)-e_0(\SIJ)=\Lambda(I,J)-e_1(I).
$$
If $\KIJ\ne 0$ then $deg(p_{\KIJ})=d-1$. Notice that
$deg(p_{\SIJ})=deg(p_{\MIJ})=d-1$, so $e_0(\KIJ)\ge 0$ and then
$\Lambda(I,J)\ge e_1(I)$.
\end{proof}

\bigskip
\begin{proposition}
\label{pasapas} \noindent (i) For all $p\ge 0$
$$
e_0(\SIJ_{[p]})=length_R(I^{p+1}/J I^p)- e_0(\KIJ_{[p]})\ge 0,
$$

\noindent and
$$
e_1(I)=\sum_{p\ge 0} (length_R(I^{p+1}/J I^p)- e_0(\KIJ_{[p]})).
$$

\noindent (ii) For all $p\ge 0$
$$
length_R(I^{p+1}\cap J/J I^p) \ge e_0(\KIJ_{[p]}),
$$

\noindent and
$$
\delta(I,J)=e_0(\KIJ)=\sum_{p\ge 0} e_0(\KIJ_{[p]}) \ge 0.
$$
\end{proposition}
\begin{proof}
$(i)$ From \propref{multiplicities} we have
$$
e_1(I)=e_0(\SIJ)=\sum_{p\ge 0}^{*} e_0(\SIJ_{[p]}),
$$
here $*$ stands for  the integers $p$ such that $deg(p_{\SIJ_{[p]}})=d-1$.

From the exact sequence of $R$-modules
$$
0 \longrightarrow \KIJ_{[p]} \longrightarrow \frac{I^{p +1}}{J
I^{p}}[T_1,\ldots,T_d] \longrightarrow \SIJ_{[p]} \longrightarrow 0,
$$
we deduce that if $deg(p_{\SIJ_{[p]}})< d-1$ then
 $length_R(I^{p+1}/J I^p)=
e_0(\KIJ_{[p]})$. Let us assume $deg(p_{\SIJ_{[p]}})=d-1$, from the additivity
of the multiplicity we deduce
$$
e_0(\SIJ_{[p]})=length_R(I^{p+1}/J I^p)- e_0(\KIJ_{[p]}),
$$
so
$$
e_1(I)=e_0(\SIJ)={\sum_{p\ge 0}} ^{*} e_0(\SIJ_{[p]}) =\sum_{p\ge 0}
(length_R(I^{p+1}/J I^p)- e_0(\KIJ_{[p]})).
$$

\noindent $(ii)$ From the well known result $gr_J(R)\cong R/J[T_1,\ldots,T_d]$
it is easy to prove that $\KIJ_{[p]}$ is in fact a submodule of
$$
\frac{I^{p +1}\cap J}{J I^{p}}[T_1,\ldots,T_d].
$$
\noindent From this we deduce
$$
length_R(I^{p+1}\cap J/J I^p) \ge e_0(\KIJ_{[p]}).
$$
 The second part of the claim follows from
\propref{multiplicities} $(iii)$
$$\delta(I,J)=\Lambda(I,J)-e_1(I)=e_0(\KIJ)=\sum_{p\ge 0}e_0(\KIJ_{[p]})\ge 0.$$
\end{proof}

\begin{remark}
\label{unificacio} Notice that some of the results on $T_{k,n}=\KIJ_{(k+n,n)}$
of \cite{Wan97d}, \cite{Wan97a}, \cite{Wan01}, \cite{Wan00} and \cite{Wan02b}
are corollaries of \propref{multiplicities} and \propref{pasapas}.
\end{remark}

Let us recall that from \cite{Wan00} and  \cite{HM95} we have
$$
\Delta(I,J) \ge \delta(I,J)=\Lambda(I,J)-e_1(I)\ge 0.
$$
In the next result we show that these inequalities can we deduced from some
"local" inequalities. For all $p\ge 0$ we define the following the integers

\medskip
\noindent $ \Delta_p(I,J)=length_R(I^{p+1}\cap J/J I^p), $

\medskip
\noindent $ \delta_p(I,J)=e_0(\KIJ_{[p]}), $ and

\medskip
\noindent $ \Lambda_p(I,J)=length_R(I^{p+1}/J I^p). $

\medskip
\noindent From the last result we deduce:

\begin{proposition}
\label{localversion} For all $p\ge 0$ the following inequalities hold
$$
\Delta_p(I,J) \ge \delta_p(I,J)=\Lambda_p(I,J)-e_0(\SIJ_{[p]})\ge 0.
$$
\noindent Summing up these inequalities with respect $p$ we get
$$
\Delta(I,J) \ge \delta(I,J)=\Lambda(I,J)-e_1(I)\ge 0.
$$
\end{proposition}

\bigskip

\section{On the depth of the blow-up algebras}

The aim of this section is to prove a refined version of  Wang's conjecture by
considering some special configurations of the set $\{\delta_p(I,J)\}_{p\ge 0}$
instead of $\delta=\sum_{p\ge 0} \delta_p(I,J)$, \thmref{generalitzacio}. As a
by-product we recover the known cases of Wang's conjecture, \corref{wang}, and
we prove a weak version of Sally's conjecture, \corref{sallyconj}.

Inspired by the proof of Polini, in the next result we generalize Claim 3 of
\cite{Pol00}.

\medskip
\begin{theorem}
\label{easykernel} Assume that $d\ge 3$. Let $I$ be an $\max-$primary ideal of
$R$ and let $J$ a minimal reduction of $I$. Let us assume that  $\KIJ \neq 0$,
and  either $\KIJ_{[p]}=0$ or
 $\KIJ_{[p]}$ is  a rank one torsion free $\res
[T_1,\ldots,T_d]$-module for $p\ge 0$. Then
$$
depth(gr_{Jt}(\rees(I)))\ge d-1.
$$
\end{theorem}
\begin{proof}
Let $p_1 < \ldots < p_n$ be the sequence of integers  such that
$\KIJ_{[p_i]}\neq 0$,  $i=1,\ldots,n$. By \lemref{vanishing} we have a finite
number of these integers. Hence, from the sequence $(\textrm{S}_{[p]})$ we get
$$
\SIJ_{[p]} \cong \MIJ_{[p]} =\frac{I^{p+1}}{JI^p}[T_1,\ldots,T_d]
$$
 for $p\neq p_1,\ldots,p_n$,
and the following exact sequences of $R-$modules
\begin{equation}
\label{aux} 0 \longrightarrow \KIJ_{[p_i]} \longrightarrow \frac{I^{p_i +1}}{J
I^{p_i}}[T_1,\ldots,T_d] \longrightarrow \SIJ_{[p_i]} \longrightarrow 0,
\end{equation}
$i=1,\ldots,n$. Notice that by hypothesis $\KIJ_{[p_i]}$ is isomorphic to an
ideal $I_i$ of $D=\res [T_1,\ldots,T_d]$, $i=1,\ldots,n$.

Let $\mathfrak p$ be a height $h \ge 2$ prime ideal of  $D$. Since
$D=\rees(J)/\max \rees(J)$, there exists a prime ideal $\mathfrak q$ of
$\rees(J)$ such that $\mathfrak p=\mathfrak q / \max \rees(J)$.

Since $depth_{\mathfrak q}(S_J(I))\ge 1$, \cite{Pol00}, by depth counting in
the set of exact sequences of $\rees(J)-$modules (SQ1) we get that
$depth_{\mathfrak q}(\SIJ_{[p_1]})\ge 1$. In fact, for $i<p_1$ we have that
$\SIJ_{[i]}\cong\MIJ_{[i]}$, so $depth_{\mathfrak q}(\SIJ_{[i]})\ge d> 2$.
Thus, if $p_1\ne 1$ then
$$
depth_{\mathfrak q}(C_2)\ge min\{depth_{\mathfrak q}(C_1), depth_{\mathfrak
q}(\SIJ_{[1]})-1\}\ge 1.
$$
 Hence, while $i< p_1$
we have that $depth_{\mathfrak q}(C_{i+1})\ge 1$, and it implies that
$$
depth_{\mathfrak q}(\SIJ_{[p_1]})\ge min\{ depth_{\mathfrak q}(C_{p_1}),
depth_{\mathfrak q}(C_{p_1+1})+1\}\ge 1.
$$
 Otherwise,
if $p_1=1$ then
$$
depth_{\mathfrak q}(\SIJ_{[1]})\ge min\{ depth_{\mathfrak q}(C_1),
depth_{\mathfrak q}(C_2)+1\}\ge 1.
$$
Hence we have $depth_{\mathfrak q}(\SIJ_{[p_1]})\ge 1$.

Depth counting on \eqref{aux} yields
$$
depth_{\mathfrak q}(\KIJ_{[p_1]})\ge min\{depth_{\mathfrak
q}(\frac{I^{p_1+1}}{JI^{p_1}}[T_1,\ldots,T_d]), depth_{\mathfrak
q}(\SIJ_{[p_1]})+1\}\ge 2
$$
because $depth_{\mathfrak q}(\frac{I^{p_1+1}}{JI^{p_1}}[T_1,\ldots,T_d])\ge d$
and $depth_{\mathfrak q}(\SIJ_{[p_1]})\ge 1$. Then
$$
depth(I_1)_{\mathfrak p}= depth(I_1)_{\mathfrak q}\ge depth_{\mathfrak q}
(\KIJ_{[p_1]})\ge 2,
$$
\cite{MatCA}. In particular we have that $\mathfrak p \notin Ass_D(I_1)$,
because $depth_{\mathfrak p}(I_1)\ge 2$, so $I_1$ is an unmixed ideal of $D$ of
height one. In fact, all the associated primes of $I_1$ have height $\le 1$,
and being  $D$ a domain and $I_1\ne 0$, the associated ideals of $I_1$ are
height $1$ ideals. Since $D$ is factorial we deduce that $I_1 \subset D$ is
principal,
 and then $depth(\KIJ_{[p_1]})=d$.

Since $depth_{\mathfrak q}(\frac{I^{p_1+1}}{JI^{p_1}}[T_1,\ldots,T_d])\ge d$
and $depth_{\mathfrak q}(\KIJ_{[p_1]})=d$, by depth counting on \eqref{aux}, we
deduce that $depth_{\mathfrak q}(\SIJ_{[p_1]})\ge d-1\ge 2$.

By depth counting on (SQ1) we get that $depth_{\mathfrak q}(C_{p_1 +1})\ge 1$.
In fact,
$$
depth_{\mathfrak q}(C_{p_1+1}) \ge min\{depth_{\mathfrak q}(C_{p_1}),
depth_{\mathfrak q}(\SIJ_{[p_1]})-1 \} \ge 1,
$$
since $ depth_{\mathfrak q}(C_{p_1})\ge 1.$ We can iterate the process and we
get that $depth_{\mathfrak q}(\SIJ_{[p]})\ge d-1$ for all $p$, in particular we
get that the $\rees(J)-$module $\SIJ$ verifies
$$depth(\SIJ_{[p]})\ge d-1.$$

\medskip
From the last row of the sequence (SQ2) we get that
$$
depth(\SIJ_{\ge r-2})\ge min\{ depth(\SIJ_{[r-2]}), depth(\SIJ_{[r-1]})\}\ge
d-1.
$$
Iterating this process in (SQ2)  we deduce that
$$
depth(\SIJ)\ge d-1.
$$

Now, let consider the exact sequence of $\rees(J)-$modules
$$
0 \longrightarrow \rees(J) \longrightarrow gr_{Jt}(\rees(I)) \longrightarrow
\SIJ  \longrightarrow 0,
$$
depth counting in this sequence give us the claim, because
$$
depth(gr_{Jt}(\rees(I)))\ge  min\{depth(\rees(J)),
 depth(\SIJ)\}
$$
 with $depth(\rees(J))=d+1$ and $depth(\SIJ)\ge
d-1$.
\end{proof}


\bigskip
\begin{lemma}
\label{torsionfree} If $\delta_p(I,J)=1$ then $\KIJ_{[p]}$ is a rank one
torsion free $\res[T_1,\ldots,T_d]$-module.
\end{lemma}

\begin{proof}
Let us recall that
$$
\KIJ_{[p]} \subset \MIJ_{[p]}=\frac{I^{p+1}}{J I^p}[T_1,\ldots,T_d]
$$
We denote by $\KIJ_{[p],j}$ the homogeneous piece of degree $j$ of $\KIJ_{[p]}$
with respect $T_1,\ldots,T_d$. Since $\frac{I^{p+1}}{J I^{p}}$ is a finite
length $R$-module there exists a composition series
$$
0=N_l \subset N_{l-1} \subset \dots \subset N_0=\frac{I^{p+1}}{J I^{p}}
$$

\noindent such that $N_i$ is a sub-$R$-module of $I^{p+1}/J I^{p}$ and
$N_i/N_{i+1} \cong \res$ for all $i=0,\dots, l-1$, i.e. $length_R(I^{p+1}/J
I^{p})=l$. Hence we have  a sequence of $R[T_1,\ldots,T_d]-$modules
$$
0=N_l[T_1,\ldots,T_d]  \subset N_{l-1}[T_1,\ldots,T_d] \subset \dots \subset
N_0[T_1,\ldots,T_d]=\frac{I^{p+1}}{J I^{p}}[T_1,\ldots,T_d].
$$
If we denote by
$$
W_{i}=\frac{\KIJ_{[p]} \cap N_i[T_1,\ldots,T_d]}{\KIJ_{[p]} \cap
  N_{i+1}[T_1,\ldots,T_d]}
\subset \frac{N_i[T_1,\ldots,T_d]}{ N_{i+1}[T_1,\ldots,T_d]}=
\res[T_1,\ldots,T_d]
$$
Since $e_0(\KIJ_{[p]})=1$ we have that $\KIJ_{[p]}\ne 0$. If $W_i=0$ for all
$i=0,\ldots, l$ then $\KIJ_{[p]}=0$, so there exists a set of indexes $0\le i_1
\le \ldots\le i_s \le l$ such that $W_{i_{j}}\neq 0$, $j=1,\ldots, s$.

We denote by $W_{i,m}$ the degree $m$ piece of $W_i$ with respect
$T_1,\ldots,T_d$. Let $m_0$ be an integer such that for all $j=1,\ldots,s$ we
have $W_{i_{j},m_0}\neq 0$. This integer exists. In fact, for each $W_{i_j}\ne
0$, there exists an integer $m_{i_j}$ such that $W_{i_j, m_{i_j}}\ne 0$. Then,
for any $t\ge 0$ we have that $W_{i_j,m_{i_j}+t}\ne 0$. So, we can choose the
maximum of these $m_{i_1},\ldots,m_{i_s}$. Then we have
\begin{eqnarray*}
length_R(\KIJ_{[p],m}) & = & \sum_{j=1}^{s} length_R (W_{i_{j},m}) =
\sum_{j=1}^{s}length_R(W_{i_{j},m_0}\cdot
W_{i_{j},m-m_0}) \\
& \ge &
\sum_{j=1}^{s}length_R(W_{i_{j},m-m_0}) \\
& \ge & s \binom{m-m_0+d-1}{d-1}
\end{eqnarray*}
for all $m\ge m_0$.

Since $e_0(\KIJ_{[p]})=1$ we deduce that $s=1$, because
$deg(p_{\KIJ_{[p]}})=d-1$ and the binomial number is a polynomial of the same
degree. Hence, $W_{i}=0$ for all $i\ne i_1$. Therefore, $\KIJ_{[p]} \cap
N_i[T_1,\ldots,T_d]=0$ for $i=i_1+1,\ldots,l$ and $\KIJ_{[p]} \cap
N_i[T_1,\ldots,T_d]=\KIJ_{[p]}$ for $i=0,\ldots, i_1$. From this we get that $
\KIJ_{[p]} = W_{i_1} \subset \res[T_1,\ldots,T_d], $ so $\KIJ_{[p]}$ is a rank
one torsion free $\res[T_1,\ldots,T_d]$-module.
\end{proof}

\bigskip

Next theorem  is the main result of this paper. We prove  a refined version of
Wang's conjecture by considering some special configurations of the set
$\{\delta_p(I,J)\}_{p\ge 0}$ instead of $\delta=\sum_{p\ge 0} \delta_p(I,J)$.
Let us consider $\bar{\delta}(I,J)$ the maximum of the integers $\delta_p(I,J)$
for $p\ge 0$.


\begin{theorem}
\label{generalitzacio} Assume that $\bar{\delta}(I,J)\le 1$. Then
$$
depth(\rees(I))\ge d-\bar{\delta}(I,J)
$$
and $depth(gr_I(R))\ge d-1-\bar{\delta}(I,J).$
\end{theorem}
\medskip
\begin{proof}
If $\bar{\delta}(I,J)=0$ then $\KIJ=0$, \propref{multiplicities}. The exact
sequence $(S)$ shows that $\MIJ \cong \SIJ$ as $A-$modules, so $depth(\SIJ)=d$,
\propref{assoc}. Let us consider the exact sequence of $A-$modules
\begin{equation*}
0 \longrightarrow \SIJ \longrightarrow gr_{Jt}(\rees(I)) \longrightarrow
\rees(J)\longrightarrow 0
\end{equation*}
Depth counting shows that $depth(gr_{Jt}(\rees(I))) \ge d$, because
$depth(\rees(J))=d+1$ and $depth(\SIJ)=d$. Hence $depth(\rees(I)) \ge d$. If
$gr_I(R)$ is a Cohen-Macaulay ring then $depth(gr_I(R))=d\ge d-1$. Otherwise,
when $depth(gr_I(R))<d=depth(R)$, from \cite{HM94b},
$$
depth(gr_I(R))=depth(\rees(I))-1\ge d-1
$$
and we get the claim.

Suppose that $\bar{\delta}(I,J)=1$. For cases $d=1,2$ we have the claim easily.
So we assume $d\ge 3$. From \lemref{torsionfree} and \thmref{easykernel} we get
$depth(gr_{Jt}(\rees(I)))\ge d-1$, so $depth(\rees(I))  \ge d-1$. Now, if
$gr_I(R)$ is Cohen-Macaulay, then $depth(gr_I(R))=d\ge d-2$. Otherwise, if
$depth(gr_I(R))< d=depth(R)$ then, \cite{HM94b},
$$
depth(gr_I(R))=depth(\rees(I))-1\ge (d-1)-1=d-2
$$
 and we get the
claim.
\end{proof}

\bigskip

\begin{remark}
The example of Wang in \cite{Wan02b}, Example 3.13, shows that the last result
is sharp in the sense that we cannot expect to have $depth(gr_{I}(R)) \ge d -1$
provided $\bar{\delta}(I,J)=1$. Precisely, this is a counterexample for the
question formulated by Guerrieri in \cite{Gue93}. She asked if it were true
that $depth(gr_I(R))\ge d-1$ for an $\max-$primary ideal $I$ in a
$d-$dimensional Cohen-Macaulay ring provided that $\Delta_p(I,J)\le 1$ $\forall
p\ge 1$. Wang reformulate the question in the regular case. Relating to this,
we are able to improve the bound for the Cohen-Macaulay case:
\end{remark}

\begin{proposition}\label{GUE-question} Let $I$
an $\max-$primary ideal and $J$ a minimal reduction of $I$ in a Cohen-Macaulay
$d-$dimensional local ring. Assuming that $\Delta_p(I,J)\le 1$ for all $p\ge 1$
we have that
$$depth(gr_I(R))\ge d-2.$$
\end{proposition}
\begin{proof}
Observe that $\delta_p(I,J)\le\Delta_p(I,J)\le 1$ for all $p\ge 1$. For $p=0$
we have that $\delta_0(I,J)=e_0(K_{[0]}^{I,J})=0$ for \lemref{vanishing} and
use \thmref{generalitzacio}.
\end{proof}

\bigskip
In the following corollary, we are able to prove the Conjecture of Wang in the
known cases, \cite{Wan00}, using the previous results and the bigraded modules
defined before. Notice that in general we have
$$
\delta(I,J) \ge \bar{\delta}(I,J),
$$
so from the last result we deduce:

\medskip
\begin{corollary}
\label{wang} Let $I$ be an $\max$-primary ideal and $J$ a minimal reduction of
$I$. Then it holds
$$
depth(gr_{I}(R)) \ge d -1 -\delta(I,J)
$$
for $\delta(I,J)=0,1$.
\end{corollary}
\medskip
\begin{proof}
Notice that $\delta(I,J)\ge \bar{\delta}(I,J)$.
\end{proof}

\bigskip
Let us recall that Valabrega and Valla  characterized under which conditions
$gr_I(R)$ is Cohen-Macaulay. They proved that given a minimal reduction  $J$ of
$I$ then $gr_I(R)$ is Cohen-Macaulay if and only if for all $n\ge 0 $ the
$n-$th Valabrega-Valla's condition holds, \cite{VV78}:
$$
 I^{n} \cap J = I^{n-1}\, J,
 $$
 i.e. $\Delta_{n-1}(I,J)=0$.

\bigskip
In the next result we prove  a weak version of Sally's conjecture,
\cite{CPP97}, \cite{Eli99}, \cite{Ros00}.

\begin{corollary}
\label{sallyconj}. Let $I$ be an $\max-$primary ideal of $R$ with minimal
reduction $J$. If $I^n \cap J= I^{n-1}J$ for $n=2,\ldots,t$ and
$length(\frac{I^{t+1} }{J I^{t}})=\epsilon \le \text{ Min} \{1, d-1\}$ then it
holds
$$
d-1- \epsilon\le depth(gr_I(R)) \le d.
$$
\end{corollary}
\begin{proof}
From \propref{localversion} we get that $\delta_{p}(I,J)=0$ for all $p\le t-1$.

Let us consider the  finitely generated $R/J$-algebra $\mathcal B$
$$
\mathcal B= \frac{\rees (I)}{J + J  t \rees (I)} =
\frac{R}{J}\oplus\bigoplus_{n \ge 1} \frac{I^{n}}{J I^{n-1}} t^{n}.
$$
Notice that $\mathcal B_{\ge 1}$ is the positive part of the degree zero piece
with respect $U$ of $\SIJ$. We can consider the Hilbert function of $\mathcal
B$, $n \ge 0 $,
$$
h_{\mathcal B}(n)= length_{R/J}(I^{n}/J I^{n-1}).
$$

From \cite{Bla95} and \cite{BN96} we deduce that
$$
h_{\mathcal B}(t+1+n)=length(\frac{I^{t+1+n} }{J I^{t+n}})\le \epsilon \le 1
$$
for all $n\ge 0$, so  $\delta_{p}(I,J)\le 1$ for all $p\ge  t$,
\propref{localversion}.

We know that $\delta_{p}(I,J)\le 1$ for all $p\ge 0$, from
\thmref{generalitzacio} we get the claim.
\end{proof}

\bigskip

\bibliographystyle{amsplain}

\providecommand{\bysame}{\leavevmode\hbox to3em{\hrulefill}\thinspace}

\end{document}
\end